\documentclass[10pt]{article}
\usepackage{amssymb}
\usepackage{amsmath}
\usepackage[dvips]{graphics}
\usepackage{epsfig}
\newtheorem{theorem}{Theorem}[section]

\newtheorem{conjecture}[theorem]{Conjecture}
\newtheorem{corollary}[theorem]{Corollary}

\newtheorem{remark}[theorem]{Remark}

\numberwithin{equation}{section}

\begin{document}
\date{}
\author{ Aristides V. Doumas \\
%EndAName
Department of Mathematics\\
National Technical University of Athens\\
Zografou Campus\\
157 80 Athens, GREECE\\
adou@math.ntua.gr}
\title{On the minimum of independent collecting processes via the Stirling numbers of the second kind}
\maketitle
\begin{abstract}
We consider the combinatorial problem where $p$ players aim to a complete set of $N$ different types of items (species) which are uniformly distributed. Let the random variables $T_{N(i)},\,\,i=1,2,\cdots,p$ denoting the number of trials needed until all $N$ types are detected (at least once), respectively for each player. This paper studies the impact of the number $p$ in the asymptotics of the expectation, the second moment, and the variance of the random variable
\begin{equation*}
M_{N(p)}: = \bigwedge_{i=1}^p T_{N(i)},\,\,\,\,\,\,N\rightarrow \infty.
\end{equation*}
The main ingredient in the expression of these quantittes are sums involving the Stirling numbers of the second kind; for which the asymptotics are explored. At the end of the paper we conjecture on a remarkable \textit{combinatorial identity}, regarding alternating binomial sums. These sums have been studied (mainly) by P. Flajolet due to their applications to digital search trees and quadtrees.
\end{abstract}

\textbf{Keywords.} Species detection; Coupon collector's problem; Stirling numbers of the second kind, digital search trees, quadtrees.\\\\
\textbf{2010 AMS Mathematics Classification.} 11B73, 34M30, 05A19.

\section{Introduction}
The \emph{coupon collector's problem (CCP)} in its classic form refers to a population whose members are of $N$ different \emph{species}. For $1 \leq j \leq N$ we denote by $p_j$ the probability that a member of the population is of type $j$, where $p_j > 0$ and $\sum_{j=1}^{N}p_{j}=1$. We refer to the $p_j$'s as the coupon probabilities. The members of the population are sampled independently \textit{with replacement} (alternatively, the population is assumed very large) and their types are recorded.

Let $T_{N(1)}, T_{N(2)},\cdots, T_{N(p)}$, the random variables denoting the number of trials needed until all $N$ types are detected (at least once) for each one of $p$ independent collectors (here and in what follows $p$ is a fixed positive integer). A main task of this note is to study the expectation of the random variable
\begin{equation*}
M_{N(p)}: = \bigwedge_{i=1}^p T_{N(i)},\,\,\,\,\,\,N\rightarrow \infty,
\end{equation*}
 i.e., the minimum of the random variables $\left\{T_{N(1)}, T_{N(2)},\cdots, T_{N(p)}\right\}$,
when the coupon probabilities are uniformly distributed, namely when
$p_{j}= 1/N$, $j=1,2,\cdots, N$ . Notice that for $p=1$ we have the classic version of the CCP. Since $M_{N(p)}$ is a non negative random variable
\begin{equation}
E\left[\,M_{N(p)}\,\right]=\sum_{k=1}^{\infty} P\left\{ M_{N(p)}\geq k\right\}.
\label{1}
\end{equation}
By independence we have
\begin{equation}
P\left\{ M_{N(p)}\geq k\right\}=P\left\{ T_{N}\geq k\right\}^{p},\qquad k=1,2,...\,,
\label{2}
\end{equation}
where $T_N$ is the random variable denoting the number of trials one collector needs until all $N$ different types are detected. Let us consider the probability $P\left\{ T_{N}\geq k\right\}$ for general values of $p_{j}$. Clearly, $k\geq N$. For each $j \in \{1,...,N\}$ it is convenient to introduce the event $A_{j}^{k}$, that the type $j$ is \underline{\textbf{not}} detected until trial $k$ (included). Then
\begin{equation*}
P\left\{ T_{N}\geq k\right\} =P\left( A_{1}^{k-1}\cup \cdot \cdot \cdot \cup
A_{N}^{k-1}\right).
\end{equation*}
By invoking the inclusion-exclusion principle one gets (see, e.g. \cite{DPM})
\begin{equation*}
P\left\{ T_{N}\geq k\right\} =\sum_{\substack{ J\subset \left\{
1,...,N\right\}  \\ J\neq \emptyset }}\left( -1\right) ^{|J|-1}\left[
1-\left( \sum_{j\in J}p_{j}\right) \right] ^{k-1},
\end{equation*}
where the sum extends over all $2^N - 1$ nonempty subsets $J$ of $\left\{
1,...,N\right\}$, while $|J|$ denotes the cardinality of $J$. Setting $p_{j}=1/N$ we get
\begin{equation}
P\left\{ T_{N}\geq k\right\} = \left(-1\right)^{N-1}\sum_{n=0}^{N-1}(-1)^{n}\binom{N}{n}\left(\frac{n}{N}\right)^{k-1}.
\label{pr1}
\end{equation}
Recall that the Stirling numbers of the second kind count the ways to partition a set of $k$ labeled objects into $N$ nonempty unlabeled subsets and they can be calculated by the so-called Euler's formula for Stirling numbers (see, e.g. \cite{GOULD}, pp.118--119):
\begin{equation*}
S(k,N)= \frac{1}{N!}\sum_{n=0}^{N}(-1)^{N-n}\binom{N}{n}n^{k}.
\end{equation*}
By invoking the Stirling numbers of the second kind in (\ref{pr1}) we  have
\begin{equation}
P\left\{ T_{N}\geq k\right\} = 1-S(k-1,N)\frac{N!}{N^{k-1}}.
\label{3}
\end{equation}
In view of (\ref{2}) and (\ref{3}), relation (\ref{1}) yields the interesting formula
\begin{equation}
E\left[\,M_{N(p)}\,\right]=\sum_{k=0}^{\infty}\left(1-S(k,N)\frac{N!}{N^{k}}\right)^{p},
\label{4}
\end{equation}
where, of course, $S(k,N)=0$ for $k<N$. The rest of our analysis is, mainly, devoted to the asymptotics of $E\left[\,M_{N(p)}\,\right]$ as $N \rightarrow \infty$.

\section{Asymptotic analysis}
Fix a positive integer $N$. Let us set
\begin{equation}
c_N:=\min\limits_{j\in \mathbb{Z^{+}}} \left\{\frac{N+j}{\ln\left(N+j\right)}>N\right\}.
\label{c}
\end{equation}
In words $c_N$ is the smallest positive integer such that the fraction above is greater than $N$. We have
\begin{equation}
E\left[\,M_{N(p)}\,\right]=\sum_{k=0}^{N+c_N-1}\left(1-S(k,N)\frac{N!}{N^{k}}\right)^{p}
+\sum_{k=N+c_N}^{\infty}\left(1-S(k,N)\frac{N!}{N^{k}}\right)^{p}.
\label{4a}
\end{equation}
The idea behind (\ref{4a}) is that the asymptotic behavior of the Stirling numbers of the second kind is known thanks to  Erd\H{o}s and Szekeres when
\begin{equation}
N<\frac{k}{\ln k}
\label{200}
\end{equation}
and thanks to \cite{L} otherwise.\footnote{For more results regarding the behavior of the Stirling numbers of the second kind, we refer the interested reader to \cite{L}, \cite{M}, and  \cite{T}.} (Notice that the inequality (\ref{200}) may also be written in terms of the Lambert $W$ function). From here and it what follows we will call the first of the sums of (\ref{4a}) as $S_{1}(N)$ and the second one as $S_{2}(N)$. We start with $S_{2}(N)$. In case where (\ref{200}) holds, we have (see \cite{S}, pp.164)
\begin{equation}
S(k,N)= \frac{N^{k}}{N!} \exp\left[\left(\frac{k}{2N}-N\right) e^{-\frac{k}{N}}\right]\left(1+o(1)\right). \label{stir}
\end{equation}
Using (\ref{stir}) in $S_{2}(N)$ (of (\ref{4a})), and from the comparison of sums and integrals we get
\begin{align}
S_{2}(N)=
\int_{N+c_N}^{\infty}\left(1-\exp\left[\left(\frac{x}{2N}-N\right) e^{-\frac{x}{N}}\right]\left(1+o(1)\right)\right)^{p}dx\,\bigg(1+o(1)\bigg).
\label{000}
\end{align}
Let us consider the integral
\begin{equation}
J(N;p):=\int_{N+c_N}^{\infty} \left(1-\exp\left[\left(\frac{t}{2N}-N\right) e^{-\frac{t}{N}}\right]\right)^{p}dt.
\label{51}
\end{equation}
Changing the variables as $e^{\frac{N+c_N}{N}-\frac{t}{N}}=x$ and integrating by parts yields
\begin{align*}
J(N;p)&=-e^{-\frac{N+c_N}{N}}\,p N\\
&\times \int_{0}^{1}\left(N+\frac{1}{2}\ln x++\frac{c_N}{2N}\right) e^{-e^{-\frac{N+c_N}{N}}\left(N+\frac{\ln x}{2}-\frac{N+c_N}{2N}\right)x}\\ &\,\,\,\,\,\,\,\,\,\,\,\,\,\,\,\,\,\,\,\,\,\,\,\,\,\,\,\,\,\,\,\,\,\,\,\,\,\,\,\,\,\,\,\,\,\,
\,\,\,\,\,\,\,\,\,\,\,\,\,\times\left(1-e^{-e^{-\frac{N+c_N}{N}}\left(N+\frac{\ln x}{2}-\frac{N+c_N}{2N}\right)x}\right)^{p-1} \ln x\,dx.
\end{align*}
Thanks to the binomial theorem we get ($p$ is always a positive integer)
\begin{align}
J(N;p)&=-e^{-\frac{N+c_N}{N}}p N\sum_{j=0}^{p-1}(-1)^{j}\binom{p-1}{j}\int_{0}^{1}\left[N\left(\ln x\right)+\frac{1}{2}\left(\ln x\right)^{2}+\frac{c_N\ln x}{2N}\right] \nonumber\\
&\times\left(1-\left(j+1\right)e^{-\frac{N+c_N}{N}}\frac{x\ln x}{2}+\left(j+1\right)e^{-\frac{N+c_N}{N}}\frac{N+c_N}{2N}x
+O\left(x\ln x\right)^{2}\right) \nonumber\\
&\times e^{-\left(j+1\right)xNe^{-\frac{N+c_N}{N}}}dx,
\label{6}
\end{align}
where we have also used the Taylor expansion of the exponential around $x=0$. Set
\begin{align*}
I_{1}(N):=&\int_{0}^{1}\left(\ln x\right) e^{-\left(j+1\right)xNe^{-\frac{N+c_N}{N}} }dx.
\end{align*}
Changing the variables as $\left(j+1\right)xNe^{-\frac{N+c_N}{N}}  =y$ yields
\begin{align*}
I_{1}(N)=\frac{e^{\frac{N+c_N}{N}}}{\left(j+1\right) N}
\times &\left[\int_{0}^{\infty} e^{-y}\ln y\,dy-\int_{\left(j+1\right)Ne^{-\frac{N+c_N}{N}} }^{\infty} e^{-y}\ln y\,dy\right.\nonumber\\
&\left.-\ln[\left(j+1\right)Ne^{-\frac{N+c_N}{N}} ]
\left(1-e^{-Ne^{-\frac{N+c_N}{N}} }\right)\right].
%\label{8}
\end{align*}
But
\begin{equation*}
\int_{0}^{\infty} e^{-y}\ln y\,dy=-\gamma,
\end{equation*}
where $\gamma=0.5772...$ is the Euler--Mascheroni constant (see, e.g. \cite{B}).
Hence,
\begin{equation}
I_{1}(N)=\frac{e^{\frac{N+c_N}{N}}}{\left(j+1\right) N}\left[-\ln N-\gamma-\ln\left(1+j\right)+1+\frac{c_N}{N}+O\left(e^{-Ne^{-\frac{N+c_N}{N}}}\ln N\right)\right].
\label{-1}
\end{equation}
%and similarly one has
%\begin{equation}
%I_{2}(N):=\int_{0}^{1}\left(\ln x\right)^{2} e^{-\left(j+1\right)e^{-1} Nx}dx=O\left(\frac{\left(\ln N\right)^{2}}{N}\right),
%\,\,N\rightarrow \infty.
%\label{110}
%\end{equation}
Working in a similar way with the rest of the integrals appearing in (\ref{6}) and finally invoking (\ref{-1})  in (\ref{6}) we get
\begin{align}
J(N;p)=p N\sum_{j=0}^{p-1}(-1)^{j}\binom{p-1}{j}\left[\frac{\ln N+\gamma-1 +\ln\left(1+j\right)}{1+j}-\frac{c_N}{N\left(j+1\right)}\right.\nonumber\\
\left.+O\left(\frac{\left(\ln N\right)^{2}}{N}\right)\right].
\label{11}
\end{align}
Next observe that
\begin{align}
&\sum_{j=0}^{p-1}(-1)^{j}\binom{p-1}{j}\frac{1}{j+1}\nonumber\\
=&\sum_{j=0}^{p-1}(-1)^{j}\binom{p-1}{j}\left(\int_{0}^{1}x^{j}dx\right)\nonumber\\
=&\int_{0}^{1}\left(1-x\right)^{p-1}dx=\frac{1}{p}.
\label{ppp1}
\end{align}
Relation (\ref{ppp1}) simplifies (\ref{11}). On the other hand let us define the constant
\begin{equation}
c_{p}:=\sum_{j=0}^{p-1}(-1)^{j}\binom{p-1}{j}\frac{\ln\left(1+j\right)}{j+1}.\label{cp}
\end{equation}
For example if $p=2$ (namely, the case we have two independent collectors), then $c_{2}= -\ln2/2$. If $p=3$ we have $c_{3}=-\ln2+(\ln3/3)$, if $p=5$ then $c_{5}=-4\ln2+2\ln3+(\ln5/5)$, etc. It is worth mentioning that for large values of $p$ sums of the type of (\ref{cp}) are of importance in applied Mathematics. We will come back to this after Theorem 2.2 below. Using (\ref{ppp1}) in (\ref{11}), and by invoking (\ref{11}) and (\ref{000}) one has as $N\rightarrow \infty$
\begin{equation}
S_{2}(N)=N\left[\ln N+\left(\gamma+p\, c_{p}-1-\frac{c_N}{N}\right)+O\left(\frac{\left(\ln N\right)^{2}}{N}\right) \right].
\label{203}
\end{equation}
Last task of our analysis is $S_{1}(N)$ of (\ref{4a}). We need the behaviour of the Stirling numbers of the second kind $S(k,N)$ when
\begin{equation}
  N\geq \frac{k}{\ln k}.
  \label{21}
\end{equation}
G. Louchard (see \cite{L}), studied this behaviour in the large deviation region, namely when
\begin{equation}
  N=k-k^{\alpha},\,\,\,a>\frac{1}{2}.
  \label{22}
\end{equation}
In particular, he proved that
\begin{equation}
  S(k,N)=\frac{1}{\sqrt{2\pi}k^{\alpha/2}}\exp\left\{  k^{\alpha}\left[\left(2-\alpha\right)\ln k+\left(1-\ln2\right)+O\left(  k^{\alpha-1}\right)\right]\right\}\left(1+O\left(  k^{\alpha-1}\right)\right).
  \label{23}
\end{equation}
Our case is covered by (\ref{22}). By invoking (\ref{23}) in $S_{1}(N)$ of (\ref{4a}), and applying Stirling's formula we get
\begin{equation}
  S_{1}(N)=\sum_{k=0}^{N+c_N-1}\left(1+o(1)\right)^{p}=N+c_N\left(1+o(1)\right),\,\,\, N\rightarrow \infty.
  \label{24}
\end{equation}
From  (\ref{4a}),  (\ref{203}) and  (\ref{24}) we arrive to our first main result, which we state in the following
\begin{theorem}
Consider the classical coupon collector's problem and $p$ independent collectors aiming to complete a set of $N$ different types of coupons, which are uniformly distributed. Let $T_{N(1)}, T_{N(2)},\cdots, T_{N(p)}$ the random variables denoting the number of trials needed until all $N$ types are detected for each one of the $p$ collectors. If we set
\begin{equation*}
M_{N(p)}: = \bigwedge_{i=1}^p T_{N(i)},
\end{equation*}
then
\begin{equation*}
E\left[\,M_{N(p)}\,\right]=\sum_{k=0}^{\infty}\left(1-S(k,N)\frac{N!}{N^{k}}\right)^{p},
\end{equation*}
where $S(k,N)$ are the Stirling numbers of the second kind. Moreover as $N\rightarrow \infty$ we have
\begin{equation*}
E\left[\,M_{N(p)}\,\right]=N\left[\ln N+\left(\gamma+p\, c_{p}\right)+O\left(\frac{\left(\ln N\right)^{2}}{N}\right) \right],
\end{equation*}
where $c_{p}$ is the constant given by the formula
\begin{equation}
c_{p}:=\sum_{j=0}^{p-1}(-1)^{j}\binom{p-1}{j}\frac{\ln\left(1+j\right)}{j+1}.\label{013}
\end{equation}
\end{theorem}

\subsection{Second moment and Variance of $M_N(p)$}
Here we will briefly present the asymptotics of the second moment for the non negative random variable $M_N(p)$ of (\ref{000}). We have
\begin{equation*}
E\left[\,M_{N(p)}^{2}\,\right]=2\sum_{k=1}^{\infty}k P\left\{ M_{N(p)}\geq k\right\}-\sum_{k=1}^{\infty} P\left\{ M_{N(p)}\geq k\right\}.
%\label{13}
\end{equation*}
Following the steps of $E\left[\,M_{N(p)}\,\right]$ we get
\begin{align}
E\left[\,M_{N(p)}^{2}\,\right]=2&\sum_{k=0}^{N+c_N-1}k\left(1-S(k,N)\frac{N!}{N^{k}}\right)^{p}
+2\sum_{k=N+c_N}^{\infty}k\left(1-S(k,N)\frac{N!}{N^{k}}\right)^{p}\nonumber\\
&-E\left[\,M_{N(p)}\,\right].
\label{14}
\end{align}
Let $S_3(N)$ and $S_4(N)$ the sums appearing in (\ref{14}) above. Using the same approximation for the Stirling numbers of the second kind, we see that the key is to obtain asymptotics (as $N \rightarrow \infty$) for the integral
\begin{equation}
L(N;p):=2\int_{N+c_N}^{\infty} t \left(1-\exp\left[\left(\frac{t}{2N}-N\right) e^{-\frac{t}{N}}\right]\right)^{p}dt,
\label{15}
\end{equation}
which is the analog of the integral $J(N;p)$ of (\ref{51}). It is now straightforward to get asymptotics for $L(N;p)$
\begin{align*}
L(N;p)&=p N^{2}\sum_{j=0}^{p-1}(-1)^{j}\binom{p-1}{j}\\
&\,\,\,\,\,\,\,\times\left[\frac{\left(\ln N\right)^{2}+2\left(\gamma+\ln\left(1+j\right)\right)\ln N}{1+j}+\frac{2\gamma\ln\left(1+j\right)-1+\frac{\pi^{2}}{6}+\gamma^{2}}{1+j}
\right.\nonumber\\
&\left.\,\,\,\,\,\,\,\,\,\,\,\,\,\,\,\,\,\,\,\,\,\,\,\,\,\,\,\,\,\,\,\,\,\,\,\,\,\,\,\,\,\,\,\,\,\,\,\,\,\,\,\,\,\,\,\,+\frac{\left(\ln\left(1+j\right)\right)^{2}}{j+1}
+O\left(\frac{\left(\ln N \right)^{3}}{N}\right)\right]
\end{align*}
as $N\rightarrow \infty$, which in turn provides asymptotics for $S_3(N)$. On the other hand if we treat $S_{4}(N)$ as we treated $S_{1}(N)$ of (\ref{24}), we have
\begin{align*}
  S_{4}(N)=&2\sum_{k=0}^{N+c_N-1}k\left(1-S(k,N)\frac{N!}{N^{k}}\right)^{p}\\
  =&\left(N+c_N\right)\left(N+c_N-1\right)\left(1+o(1)\right)^{p}=N^2\left(1+o(1)\right),\,\,\, N\rightarrow \infty
 % \label{24}
\end{align*}
and finally, arrive at the following
\begin{theorem}
Let $M_{N(p)}$ as defined in Theorem 2.1. Then
\begin{equation*}
E\left[\,M_{N(p)}^{2}\,\right]=2\sum_{k=0}^{\infty}k\left(1-S(k,N)\frac{N!}{N^{k}}\right)^{p}-E\left[\,M_{N(p)}\,\right],
\end{equation*}
where $S(k,N)$ are the Stirling numbers of the second kind. In particular as $N\rightarrow \infty$ we have
\begin{align*}
E\left[\,M_{N(p)}^{2}\,\right]=N^{2}\left[\left(\ln N\right)^{2}+2\left(\gamma+p\, c_{p}\right)\ln N+\gamma^{2}+\frac{\pi^{2}}{6}+2 p c_{p}\gamma\right.\\
\left. + p\, w_{p}+O\left(\frac{\left(\ln N \right)^{3}}{N}\right) \right],
\end{align*}
where $c_{p}$ is the constant given in Theorem 2.1, and $w_{p}$ is the constant which may be calculated explicitly by the formula
\begin{equation}
w_{p}:=\sum_{j=0}^{p-1}(-1)^{j}\binom{p-1}{j}\frac{\left(\ln\left(1+j\right)\right)^{2}}{j+1}.\label{014}
\end{equation}
Moreover, for the variance of the random variable $M_N$ we have as $N\rightarrow \infty$
\begin{equation}
V\left[\,M_{N(p)}\,\right]\sim \left(\frac{\pi^{2}}{6}+p w_{p}-p^{2}c^{2}_{p}\right)N^{2}.
\label{16}
\end{equation}
\end{theorem}
\begin{remark}
From Theorems 2.1 and 2.2 we see that the leading term of the first and the second moment of the random variable $M_{N(p)}$ is the same with the classic version and independent of $p$, which first appears in the second term. However, $p$ \textit{appears} in the leading term of the variance. We remind the reader that
\begin{align*}
E\left[\,T_{N(1)}\,\right]&=N\sum_{j=1}^{N}\frac{1}{j}=N\left(\ln N+\gamma+\frac{1}{2N}+O\left(\frac{1}{N^2}\right)\right)\\
E\left[\,T_{N(1)}^{2}\,\right]&=N^{2}\left[\left(\sum_{j=1}^{N}\frac{1}{j}\right)^{2}
+\sum_{j=1}^{N}\frac{1}{j^2}\right]\\
&=N^{2}\left[\ln^{2} N+2\gamma\ln N+\gamma^{2}+\frac{\pi^{2}}{6}+O\left(\frac{\ln N}{N}\right)\right]\\
V\left[\,T_{N(1)}\,\right]&\sim\frac{\pi^{2}}{6}N^{2},\,\,\,\,\,\,\,\,\,\,\,\,\,N\rightarrow\infty,
\end{align*}
see, e.g., \cite{DPM}.
\end{remark}

\subsection{A few words for the case when $p$ becomes infinitely large}
The delicate problem of estimating asymptotically high order differences of some fixed
numerical sequence $\left\{f_k\right\}$
\begin{equation*}
D_n\left[f\right]:=\sum_{k=0}^{n}\binom{n}{k}(-1)^{k}f_k
\end{equation*}
goes back in mid $1960s$ to De Bruijn, Knuth, and Rice who showed their central role in the evaluation of data structures based on a binary representation of data. Applications of these binomial sums refer mainly to digital search tress (an alternative of \textit{Rice integrals}, see \cite{FFF}) and quadtrees (see, \cite{F} -- \cite{FF}).\\
Using the identity
\begin{equation*}
 \frac{1}{j+1}\binom{n-1}{j}=\frac{1}{n}\binom{n}{j+1}
\end{equation*}
the quantities $c_p$ and $w_p$ of (\ref{013}) and (\ref{014}) become
\begin{equation}
c_{p}=-\frac{1}{p}\sum_{k=1}^{p}(-1)^{k}\binom{p}{k}\ln k,\,\,\,\,\,w_{p}=-\frac{1}{p}\sum_{k=1}^{p}\binom{p}{k}(-1)^{k}\left(\ln k\right)^{2}
\label{00}
\end{equation}
respectively. By exploiting the techniques presented by P. Flajolet and R. Sedgewick (see \cite{F}) one has as $n \rightarrow \infty$:
\begin{align*}
\sum_{k=1}^{n}\binom{n}{k}(-1)^{k}\ln k=&\ln\left(\ln n\right)+\gamma+\frac{\gamma}{\ln n}-
\frac{\gamma^{2}+\frac{\pi^{2}}{6}}{2\left(\ln n\right)^{2}}
+O\left(\frac{1}{\left( \ln n\right)^{3}}\right)\\
\sum_{k=1}^{n}\binom{n}{k}(-1)^{k}\left(\ln k\right)^{2}=&-\left(\ln\left(\ln n\right)\right)^{2}-2\gamma\ln\left(\ln n\right)+\frac{\pi^{2}}{6}-\gamma^{2}-2\gamma\frac{\ln\left(\ln n\right)}{\ln n}\\
&+\frac{\left(\gamma^{2}+\frac{\pi^{2}}{6}\right)\ln\left(\ln n\right)}{\left(\ln n\right)^{2}}-\frac{2\gamma^{2}}{\ln n}+\frac{\gamma^{2}-\frac{\pi^{2}}{6}}{\left(\ln n\right)^{2}}+O\left(\frac{\ln\left(\ln n\right)}{\left(\ln n\right)^{3}}\right).
\end{align*}
By invoking the above asymptotics in (\ref{00}) we get the very interesting result
\begin{corollary}
\begin{equation}
\lim_{p\rightarrow \infty}\left(p^{2}c^{2}_{p}-p w_{p}\right)=\frac{\pi^{2}}{6}.
\label{18}
\end{equation}
\end{corollary}
It is remarkable that $\pi^{2}/6$ appears in this limit above, which is a difference of binomial alternating sums involving logarithms. Now Theorem 2.2. implies that as the number of the independent collectors goes to infinity the variance $V\left[\,M_{N(p)}\,\right]$, naturally, vanishes.
\begin{conjecture} The sequence
\begin{equation}
a_p := \frac{\pi^{2}}{6}+p w_{p}-p^{2}c^{2}_{p}
\label{18}
\end{equation}
is decreasing in $p$.
\end{conjecture}
\begin{remark} From the conjecture above one can prove the following remarkable identity
\begin{equation}
\frac{\pi^{2}}{6}+p w_{p}-p^{2}c^{2}_{p}>0 , \,\,\,p\in \mathbb{N}.
\label{17}
\end{equation}
\end{remark}
We continue with the following\\\\
\textbf{Examples.}\\
$(i)$ The case $p=1$, namely the case of the classic coupon collector's problem. Then $c_{1}=w_{1}=0$ and Theorems 2.1--2.2 yield
\begin{align*}
E\left[\,M_{N(1)}\,\right]&=N\left[\ln N+\gamma+O\left(\frac{\left(\ln N\right)^{2}}{N}\right) \right]\\
E\left[\,M_{N(1)}^{2}\,\right]&=N^{2}\left[\left(\ln N\right)^{2}+2\gamma\ln N+\gamma^{2}+\frac{\pi^{2}}{6}+O\left(\frac{\left(\ln N\right)^{3}}{N}\right) \right]\\
V\left[\,M_{N(1)}\,\right]&\sim \frac{\pi^{2}}{6}N^{2},
\end{align*}
in accordance with Remark 2.3.\\
$(ii)$ If $p=2$, then we have two independent collectors. Hence, $c_{2}=-\ln2/2$ and $w_{2}=-(\ln2)^{2}/2$, and Theorems 2.1--2.2 immediately imply
\begin{align*}
E\left[\,M_{N(2)}\,\right]&=N\left[\ln N+\left(\gamma-\ln2\right)+O\left(\frac{\left(\ln N\right)^{2}}{N}\right) \right]\\
E\left[\,M_{N(2)}^{2}\,\right]&=N^{2}\left[\left(\ln N\right)^{2}+2\left(\gamma-\ln2\right)\ln N+\gamma^{2}+\frac{\pi^{2}}{6}-2\gamma\ln 2 -\left(\ln 2\right)^{2}\right.\\
&\left.\,\,\,\,\,\,\,\,\,\,\,\,\,\,\,\,\,\,\,\,
\,\,\,\,\,\,\,\,\,\,\,\,\,\,\,\,\,\,\,\,
\,\,\,\,\,\,\,\,\,\,\,\,\,\,\,\,
\,\,\,\,\,\,\,\,\,\,\,\,\,\,\,\,\,\,\,\,\,\,\,\,
\,\,\,\,\,\,\,\,\,\,\,\,\,\,\,\,\,\,\,\,\,\,\,\,\,\,\,\,\,\,\,\,\,\,\,\,\,\,\,\,
\,\,\,\,\,\,\,\,+O\left(\frac{\left(\ln N\right)^{3}}{N}\right) \right]\\
V\left[\,M_{N(2)}\,\right]&\sim \left(\frac{\pi^{2}}{6}-2\left(\ln 2\right)^{2}\right)N^{2}.
\end{align*}
Notice that the variance decreases when the collectors become two instead of one (as expected). We will give closure with the following
\begin{remark} Let us consider the case where $p=2$. In view of (\ref{2}) and (\ref{pr1}) relation (\ref{1}) yields
\begin{align}
&E\left[\,M_{N(2)}\,\right]=\sum_{k=1}^{\infty}
\left[ \sum_{n=0}^{N-1}(-1)^{n}\binom{N}{n}\left(\frac{n}{N}\right)^{k-1}\right]^{2}\nonumber\\
&=\sum_{n=0}^{N-1}\binom{N}{n}^{2}\left(\frac{n}{N}\right)^{2k-2}
+2\sum_{0\leq n_{1}<n_{2}\leq N-1}\left(-1\right)^{n_{1}+n_{2}}\binom{N}{n_{1}}\binom{N}{n_{2}}\left(\frac{n_{1}n_{2}}{N^{2}}\right)^{k-1}\label{19}\\
&=\sum_{n=0}^{N-1}\binom{N}{n}^{2}\frac{1}{1-\left(\frac{n}{N}\right)^{2}}
+2\sum_{0\leq n_{1}<n_{2}\leq N-1}\left(-1\right)^{n_{1}+n_{2}}\binom{N}{n_{1}}\binom{N}{n_{2}}\frac{1}{1-\frac{n_{1}n_{2}}{N^{2}}},\label{20}
\end{align}
where the last equation follows by summing the geometric series of (\ref{19}). An aside result of Theorem 2.1 is that it provides asymptotics for expressions similar to ((\ref{19})--(\ref{20})).
\end{remark}
\bigskip
\textbf{Acknowledgements.} The author is indebted to Professor Vassilis G. Papanicolaou for placing the problem and for helpful comments and constructive remarks.

\end{document}